\documentclass[reqno,onecolumn,oneside]{paper}
\usepackage[english]{babel}
\usepackage{amsmath,amsthm,amsfonts,amssymb,latexsym,mathrsfs}
\usepackage[all]{xy}
\usepackage[sort]{cite}

\newtheorem{theorem}{Theorem}[section]
\newtheorem{proposition}[theorem]{Proposition}
\newtheorem{lemma}[theorem]{Lemma}
\newtheorem{corollary}[theorem]{Corollary}
\theoremstyle{definition}
\newtheorem{definition}[theorem]{Definition}
\newtheorem{exm}[theorem]{Example}

\theoremstyle{remark}
\newtheorem{remark}[theorem]{Remark}

\begin{document}

\title{The Logic of CMV-Algebras}
\author{Antonio Di Nola, Ciro Russo \\ \footnotesize{Dipartimento di Matematica ed Informatica -- Universit\`a di Salerno, Italy} \\ \footnotesize{\texttt{$\{$adinola, cirusso$\}$@unisa.it}} \and Brunella Gerla  \\ \footnotesize{Dipartimento di Informatica e Comunicazione -- Universit\`a dell'Insubria, Italy} \\ \footnotesize{\texttt{brunella.gerla@uninsubria.it}}}
\date{}
\maketitle

\begin{center}
{\bf Dedicated to the memory of Vittorio Cafagna.}
\end{center}

\begin{abstract}
In \cite{diflge} the authors defined the CMV-algebras, an extension of MV-algebras obtained by adding a binary operation whose behaviour, with respect to the MV-algebra structure, reflects the one of composition of McNaughton functions in the single generated free MV-algebra.

In this paper, once recalled some properties of CMV-algebras, we introduce an expansion of the one-variable fragment of \L ukasiewicz propositional logic whose algebraic semantics is the variety of CMV-algebras.
\end{abstract}

\section*{Introduction}
\L ukasiewicz logic is a non-classical logic whose synctactic apparatus was introduced by \L ukasiewicz and Tarski in \cite{luktar}. A first, synctactic, proof of the completeness of \L ukasiewicz's infinite-valued sentential calculus was given by Rose and Rosser in \cite{ros}. Later on~---~in \cite{chang,chang2}~---~Chang introduced the class of MV-algebras and showed that the Lindenbaum algebra of \L ukasiewicz logic is an MV-algebra. Since then, MV-algebras were deeply investigated by many authors and they received a great impulse from the work \cite{mun}, by Mundici, where a categorical equivalence between MV-algebras and lattice-ordered Abelian groups with strong unit was established.

We recall that an MV-algebra is a structure $\mathbf{A} = \langle A, \oplus , ^*, 0\rangle$ of type (2,1,0) satisfying the equations
\begin{enumerate}
\item[-]$(x\oplus y) \oplus z = x \oplus (y\oplus z)$,
\item[-]$x \oplus y = y \oplus x$,
\item[-]$x \oplus 0 = x$,
\item[-]$(x^*)^* = x$,
\item[-]$x \oplus 0^*=0^*$,
\item[-]$(x^* \oplus y)^* \oplus y = (y^* \oplus x)^* \oplus x$.
\end{enumerate}
The element $0^*$ is denoted by $1$. The prototypical example of MV-algebra is the real unit interval $[0,1]$ equipped with the truncated sum $x \oplus y := \min\{1, x+y\}$ and the involution $x^* := 1 - x$. It is well-known that the variety of MV-algebras is generated by this algebra and that the free MV-algebra over $n$ generators is the algebra of \emph{McNaughton functions}~---~i.e. piecewise linear functions with integer coefficients~---~from $[0,1]^n$ to $[0,1]$, equipped with the operations naturally inherited from the MV-algebra over $[0,1]$.

Since the free MV-algebra is an algebra of functions, it is naturally equipped with the further operation of map composition. In \cite{diflge} the authors defined an extension of MV-algebras obtained by adding another binary operation whose behaviour with respect to the MV-algebraic structure mirrors the composition of McNaughton functions in the free algebras: the Composition MV-algebras, \emph{CMV-algebras} for short. Many results that we present in this paper can also be found in \cite{diflge}.

The introduction of this class was stimulated by the observation of a connection between the behaviour of McNaughton functions and chaotic deterministic dynamical systems. In this paper we recall some properties of CMV-algebras and introduce an expansion of the one-variable fragment of \L ukasiewicz propositional logic whose algebraic semantics is the variety of CMV-algebras.

The paper is organized as follows. In the first section we recall the definition and basic properties of CMV-algebras, showing many examples of such structures. Among them, an important role is played by the free MV-algebra in one generator, $\mathbf{M}_1$, endowed with a monoid operation that is essentially map composition and that is also the example of CMV-algebra that inspired the aforementioned paper \cite{diflge}.

In Section~\ref{cmvend} we show the link between composition MV-algebras and MV-algebra endomorphisms: the monoid reduct of any CMV-algebra $\mathbf{A}$ is isomorphic to a submonoid of the endomorphism monoid of the MV-algebra reduct of $\mathbf{A}$.

In Section~\ref{cmvcong} we treat CMV-ideals and congruences, presenting several examples of them. We also present simple algebras and point out a wide class of such algebras.

Modules over CMV-algebras are introduced in Section~\ref{mod}. Apart from many examples presented, we will see that every MV-algebra is canonically a module over the CMV-algebra $\mathbf{M}_1$.

Finally, in Section~\ref{logic}, we present the logic $\operatorname{S}\operatorname{\textnormal{\L}}_\omega^1$~---~that is an expansion of the one-variable fragment of \L ukasiewicz propositional calculus~---~whose associated Tarski-Lindenbaum algebra is a CMV-algebra.

\section{Composition MV-algebras}

In this section we introduce the definition of CMV-algebras as a special case of composition algebras.

\begin{definition}\label{compop}
Let $f: A^n \longrightarrow A$ be a n-ary operation on a set $A$ and suppose that a monoid structure $\langle A, \diamond, i \rangle$ is defined on the same set. We say that $f$ and the monoid structure are compatible if $f(x_1, x_2, \ldots, x_n) \diamond y = f(x_1 \diamond y, x_2 \diamond y, \ldots, x_n \diamond y)$. In the special case $n = 0$ we have an element $a \in A$ and the condition reduces to $a \diamond y = a$.
\end{definition}

\begin{definition}\label{compalg}
Let $A$ be an algebra of a given signature and $\langle A, \diamond, i\rangle$ a monoid on the same universe. If all the operations of $A$ are compatible with the monoid structure we say that $A$ is a \emph{composition algebra} (of the corresponding type) with respect to the given monoid.

A \emph{morphism} of composition algebras (of the same type) is a morphism of the same algebras that is also a morphism of the associated compatible monoid structures.
\end{definition}
Clearly a composition algebra is an algebra itself but the previous definition is more useful for the scope of the present work.

\begin{definition}\label{comprel}
Let $\Phi \subseteq A^n$ be an $n$-ary relation on $A$ and $\langle A, \diamond, i\rangle$ a monoid. We say that $\Phi$ and the monoid structure are compatible if from $(x_1, x_2, \ldots, x_n) \in \Phi$ follows $(x_1 \diamond y, x_2 \diamond y, \ldots, x_n \diamond y) \in \Phi$, for every $y \in A$.
\end{definition}

An important example of composition algebras is the following. Consider a group $\langle G,+\rangle$ and the set of all the functions $f: G \longrightarrow G$ equipped with the operation $+$ defined pointwise. The composition of functions $\circ$ makes the algebra $\langle G^G,+ \rangle$ a composition algebra. Note that the composition algebra $\langle G^G, +, \circ\rangle$ is a near-ring (see \cite{pilz}).

We already recalled the definition of MV-algebra in the introduction. MV-algebra morphisms are defined in the usual way according to the signature, namely they are maps that preserve $\oplus$, $^*$ and $0$; the same holds for congruences, i.e. a congruence over an MV-algebra is an equivalence relation that preserves the operations. Let us recall also some basic facts about MV-algebras.

First of all we observe that Boolean algebras are precisely those MV-algebras satisfying the additional equation $x \oplus x = x$. Every MV-algebra contains as a subalgebra the two-element boolean algebra $\{0,1\}$; the set $B(\mathbf{A})$ of all idempotent elements of an MV-algebra $\mathbf{\mathbf{A}}$ is the largest boolean algebra contained in $A$ and is called the \emph{Boolean skeleton} of $\mathbf{A}$.

On every MV-algebra $\mathbf{A}$ we define also further (derived) operations and a structure of bounded lattice. In particular, $\odot$ is the operation defined by $x \odot y = (x^* \oplus y^*)^*$, $1$ (i.e. $0^*$) is the neutral element for $\odot$ and $\to$ is defined by $x \to y = x^* \oplus y$; moreover we define $x \ominus y := x \odot y^*$ and denote by $nx$ the element $\underbrace{x \oplus \ldots \oplus x}_{\textrm{$n$-times}}$ and by $x^n$ the element $\underbrace{x \odot \ldots \odot x}_{\textrm{$n$-times}}$. The order relation on $\mathbf{A}$ is defined by $x \leq y$ if and only if $x^* \oplus y = 1$. In this way we obtain a bounded lattice (with 0 as the bottom element and 1 as the top), with $x \vee y := x \oplus (x^* \odot y)$ and $x \wedge y := x \odot (x^* \oplus y)$. Last, the structure $\langle A, \vee, \wedge, \odot, 0, 1\rangle$ is a bounded integral commutative residuated lattice and the residuum of $\odot$ is $\to$.

\begin{exm}\label{exampleMV}
\begin{enumerate}
\item The real unit interval $[0,1]$, equipped with operations
$$x \oplus y = \min\{1,x+y\}, \quad x \odot y = \max\{0,x+y-1\}, \quad x^* = 1 - x,$$
is an MV-algebra, often referred to as the \emph{standard MV-algebra}; in this algebra the order relation is the one among real numbers. For any $n \in \mathbb N$, the set $\operatorname{\textnormal{\L}}_{n+1} = \left\{0, \frac{1}{n}, \ldots, \frac{n-1}{n}, 1\right\}$, equipped with the same operations, is a finite linearly ordered MV-algebra.
\item If $X$ is any set and $\mathbf{A}$ is an MV-algebra, $\mathbf{A}^X = \langle A^X, \oplus, ^*, 0\rangle$ is an MV-algebra with the operations defined pointwisely from those in $\mathbf{A}$.
\item The set $M_n$ of all the functions from $[0,1]^n$ to $[0,1]$ that are continuous, piecewise linear and such that each linear piece has integer coefficients~---~with the operations defined, again, pointwisely from the ones in $[0,1]$~---~is an MV-algebra. Such functions are called \emph{McNaughton functions} and the MV-algebra of McNaughton functions with $n$ variables is known to be the free MV-algebra over $n$ generators.
\end{enumerate}
\end{exm}

\begin{definition}\label{mvid}
Let $\mathbf{A}$ be an MV-algebra. An \emph{MV-ideal} of $\mathbf{A}$ is a non-empty subset $I$ of $A$ such that
\begin{itemize}
\item[-] if $x \in I$ and $x \geq y$ then $y \in I$,
\item[-] if $x, y \in I$ then $x \oplus y \in I$.
\end{itemize}

Let $\mathbf{A}$ be an MV-algebra and $I \subseteq A$ an MV-ideal of $\mathbf{A}$. The equivalence relation $\sim_I$ defined by $x \sim_I y$ if and only if $(x \odot y^*) \oplus (y \odot x^*) \in I$, is a congruence of the MV-algebra $\mathbf{A}$, called the \emph{congruence generated by $I$}. Conversely, if $\sim$ is a congruence of $\mathbf{A}$ then the set $I_\sim = \{x \in A \mid x \sim 0\}$ is an MV-ideal of $\mathbf{A}$.
\end{definition}

An ideal $I$ is called \emph{prime} if it is proper, that is $I \neq A$ and, for all $x, y \in A$, either $x \ominus y \in I$ or $y \ominus x \in
I$; $I$ is called \emph{maximal} if it is proper and no proper ideal of $\mathbf{A}$ strictly contains $I$. The radical $\operatorname{Rad}(\mathbf{A})$ of $\mathbf{A}$ is the intersection of all maximal ideals of $\mathbf{A}$ and consists of $0$ and all the infinitesimals of the algebra, where an element $a$ is called an \emph{infinitesimal} if $a \neq 0$ and $na \leq a^*$ for all $n \in \mathbb N$. An MV-algebra $\mathbf{A}$ is called \emph{perfect} if $A = \operatorname{Rad}(\mathbf{A}) \cup \operatorname{Rad}(\mathbf{A})^*$, where $\operatorname{Rad}(\mathbf{A})^* := \{x \mid x^* \in \operatorname{Rad}(\mathbf{A})\}$. An MV-algebra $\mathbf{A}$ is \emph{simple} if it has only one proper ideal.

The following representation theorem holds:
\begin{theorem}\label{chang}
An equation holds in every MV-algebra if and only if it holds in the standard MV-algebra $[0,1]$.
\end{theorem}

An important result by Mundici \cite{mun} establishes an equivalence between the category of MV-algebras and the one of lattice-ordered abelian groups ($\ell$-groups) with strong unit. A standard reference for MV-algebras is \cite{cigdotmun}.

\begin{definition}\label{cmv}
A \emph{composition MV-algebra}~---~\emph{CMV-algebra} for short~---~is a structure $\mathbf{A} = \langle A, \oplus, *, 0, \diamond, i\rangle$ such that $\langle A, \oplus, ^*, 0\rangle$ is an MV-algebra, $\langle A, \diamond, i\rangle$ is a monoid and these two structures are compatible in the sense of Definition~\ref{compalg}, i.e. for every $x, y, z \in A$, the following conditions hold:
\begin{enumerate}
\item[$(i)$] $(y \oplus z) \diamond x = (y \diamond x) \oplus (z \diamond x)$,
\item[$(ii)$] $x^* \diamond y = (x \diamond y)^*$,
\item[$(iii)$] $0 \diamond x = 0$.
\end{enumerate}
\end{definition}
Observe that $(ii)$ and $(iii)$ yield also $1 \diamond x = 1$.

Since CMV-algebras are an equational class, they form a variety that we denote by $\mathcal{CMV}$. We shall denote by $\mathbf{A}_{\mathcal{MV}}$ the MV-reduct of the CMV-algebra $\mathbf{A}$, and by $\mathbf{A}_\mathcal{M}$ the monoid reduct of $\mathbf{A}$. Further, when there is no danger of confusion, we shall denote a structure by its domain.

\begin{remark}
If $0=1$ the algebra is reduced to $\{0\}$; we call that algebra the \emph{trivial} CMV-algebra. In a non trivial CMV-algebra the operation $\diamond$ is not commutative: $1 \diamond 0 = 1$ but $0 \diamond 1 = 0$.
\end{remark}

\begin{lemma}\emph{\cite{diflge}}\label{lemma:proprts}
Let $\mathbf{A}$ be a non trivial CMV-algebra. Then the following properties hold:
\begin{enumerate}
\item[$(i)$] $0 < i < 1$;
\item[$(ii)$] $i^* \diamond i^* = i$;
\item[$(iii)$] $i^* \neq i$;
\item[$(iv)$] $a \diamond x = b \diamond x$ \ for all $x$ \ if and only if \ $a = b$.
\end{enumerate}
\end{lemma}

\begin{corollary}
Every non trivial CMV-algebra has at least four elements.
\end{corollary}

\begin{proposition}\emph{\cite{diflge}}\label{prop:compa}
Let $\mathbf{A}$ be a CMV-algebra. The operations $\odot, \vee, \wedge, \ominus$ are compatible with the monoid structure. Moreover the order relation is compatible with the monoid structure.
\end{proposition}

\begin{theorem}
There are no totally ordered non trivial CMV-algebras.
\end{theorem}
\begin{proof}
If $i \leq i^*$, using Proposition~\ref{prop:compa}, we get $i \diamond i^* \leq i^* \diamond i^*$ so $i^* \leq i$ and so $i = i^*$ which is absurd by Lemma~\ref{lemma:proprts}. If $i^* \leq i$ then $i^* \diamond i^* \leq i \diamond i^*$ and we get $i \leq i^*$ so, again, $i = i^*$.
\end{proof}

\begin{corollary}
If $\mathbf{A}$ is a CMV-algebra, then the MV-reduct $\mathbf{A}_{\mathcal{MV}}$ of $\mathbf{A}$ is neither simple nor perfect.
\end{corollary}
\begin{proof}
It is known that simple MV-algebras are totally ordered; thus $\mathbf{A}_{\mathcal{MV}}$ cannot be simple. On the other hand, it is self-evident that $i, i^* \notin \operatorname{Rad}(\mathbf{A})$ for any CMV-algebra $\mathbf{A}$, hence $\mathbf{A}$ cannot be perfect.
\end{proof}

Now we come to some examples of CMV-algebras.
\begin{exm}\label{ex:CMV}
\begin{enumerate}
\item Let $\mathbf{M}$ be any MV-algebra and let $\langle M^M, \circ, i\rangle$ be the monoid of functions from $M$ to itself, with the operation $\circ$ of composition and the identity function $i$. $\mathbf{M}^M = \langle M^M, \oplus, *, 0, 1\rangle$ is an MV-algebra with pointwise defined operations. It is easy to check that $\mathbf{A} = \langle M^M, \oplus, *, 0, \circ, i\rangle$ is a CMV-algebra.
\item Let $\mathbf B_2 = \{0,1\}$ be the Boolean algebra with two elements. We refer to the previous example and set $\mathbf{M} = \mathbf B_2$ and denote by $\mathbf{A}_4$ the corresponding CMV-algebra. It is easy to show that $\mathbf{A}_4$ is, up to isomorphisms, the unique CMV-algebra with four elements.
\item Let $\operatorname{\textnormal{\L}}_3 = \{0, 1/2, 1\}$ be the standard MV-algebra with three elements and $A$ be the set of functions $f: \operatorname{\textnormal{\L}}_3 \longrightarrow \operatorname{\textnormal{\L}}_3$ such that $f[\{0,1\}] \subseteq \{0,1\}$. With pointwise operations and map composition, $A$ becomes a composition MV-algebra. Indeed it is a subalgebra of a CMV-algebra like in Example 1. In general if $\mathbf{A} \subseteq \mathbf{B} \subseteq \mathbf{C}$ are MV-algebras, then the set $\{f \mid f(x) \in A \textrm{ for all } x \in B\}$ is the domain of a subalgebra of $\mathbf{C}^C$.
\item Consider the MV-algebra $\mathbf{M}_1$ of the McNaughton functions of one variable. The composition of two McNaughton functions (in one variable) is still a McNaughton function, so $\mathbf{M}^\mathcal C_1 = \langle M_1, \oplus, *, 0, \circ, i\rangle$ is a CMV-algebra, where $i$ is the identity map of $[0,1]$.
\item The set $\Lambda$ of all the continuous maps $f: [0,1] \longrightarrow [0,1]$ can be endowed with a CMV-algebra structure as in the previous example.
\end{enumerate}
\end{exm}
An example of CMV-algebra of functions with more than one variable is given after Proposition~\ref{prop:MVtoCMV}.

\begin{definition}
A CMV-algebra $\mathbf{A}$ is called \emph{functional} if it is a subalgebra of a CMV-algebra $\mathbf{M}^M = \langle M^M, \oplus, *, 0, \circ, i\rangle$ where $\mathbf{M}$ is an MV-algebra, $\langle M^M, \oplus, *, 0, 1\rangle$ is the MV-algebra obtained by $M$ with pointwise operations and $\circ$ is the composition of functions (see Example~\ref{ex:CMV}.1).
\end{definition}

\begin{theorem}[Cayley-type theorem]
All CMV-algebras are functional.
\end{theorem}
\begin{proof}
Let $\mathbf{A} = \langle A, \oplus, *, 0, \diamond, i \rangle$ be a CMV-algebra, $\mathbf{A}_{\mathcal{MV}}$ its MV-reduct and let $\mathbf{A}^A = \langle A^A, \oplus, *, 0, \circ, i\rangle$ be the CMV-algebra of all the functions from $A$ to $A$.

Consider, for any $a \in A$, the function $f_a \in A^A$ defined by $f_a(x) = a \diamond x$ for every $x \in A$. Then the map $\mu: A \longrightarrow A^A$ defined by $\mu \left(a\right) = f_a$ is a CMV-monomorphism of $\mathbf{A}$ onto  $\mathbf{A}^A$. Hence $\mathbf{A}$ is functional.
\end{proof}

There is another way of considering the relation between an MV-algebra $\mathbf{M}$ and the CMV-algebra $\mathbf{M}^M$. We can identify the elements of $M$ with the constant functions of $M^M$. This means to consider the map $\tau: a \in M \longmapsto c_a \in M^M$, where $c_a$ are constant functions of $M$, i.e. $c_a(x) = a$ for all $x \in M$. It is easy to show that $\tau$ is an MV-monomorphism and so $\mathbf{M}$ is isomorphic to the MV-subalgebra $\tau[\mathbf{M}]$ of $\mathbf{M}^M_{\mathcal{MV}}$ consisting of constant functions. If $\mathbf{M}$ is not trivial, i.e. $0 \neq 1$, then $\tau[\mathbf{M}]$ never contains the identity map of $M^M$, so it will never be a CMV-subalgebra of $\mathbf{M}^M$, although it is closed under composition (indeed $\tau[\mathbf{M}]$ is a semigroup of ``left zeros'' with respect to composition). We can consider the CMV-subalgebra $\widetilde{\mathbf{M}}$ generated by $\tau[\mathbf{M}]$ in $\mathbf{M}^M$. In this context we have

\begin{proposition}\label{tilde}
$\widetilde{\operatorname{\textnormal{\L}}_{n+1}} = \operatorname{\textnormal{\L}}_{n+1}^{\operatorname{\textnormal{\L}}_{n+1}}$.
\end{proposition}
\begin{proof}
Since $\widetilde{\operatorname{\textnormal{\L}}_{n+1}}$ is a subalgebra of $\operatorname{\textnormal{\L}}_{n+1}^{\operatorname{\textnormal{\L}}_{n+1}}$, by definition, we just need to prove that $\operatorname{\textnormal{\L}}_{n+1}^{\operatorname{\textnormal{\L}}_{n+1}} \subseteq \widetilde{\operatorname{\textnormal{\L}}_{n+1}}$. Let us denote, for simplicity, $\operatorname{\textnormal{\L}}_{n+1} = \{0, 1, \ldots, n\}$ and a function $f: \operatorname{\textnormal{\L}}_{n+1} \longrightarrow \operatorname{\textnormal{\L}}_{n+1}$ by $(f(0), f(1), \ldots, f(n))$. We obviously have that the functions $i = (0, 1, \ldots, n)$, $i^* = (n, n-1, \ldots, 0)$, $f_c = (c, \ldots, c)$ belong to $\widetilde{\operatorname{\textnormal{\L}}_{n+1}}$ for all $c \in \operatorname{\textnormal{\L}}_{n+1}$.

Then
$$(n, n-1, \ldots, 0) \ominus (n-1, n-1, \ldots, n-1) = (1, 0, \ldots, 0) \in \widetilde{\operatorname{\textnormal{\L}}_{n+1}}$$
and
$$(0, 1, \ldots, n) \ominus (n-1, n-1, \ldots, n-1) = (0, \ldots, 0, 1) \in \widetilde{\operatorname{\textnormal{\L}}_{n+1}}.$$
So $(n, 0, \ldots, 0) \in \widetilde{\operatorname{\textnormal{\L}}_{n+1}}$ and, from
$$(n, n-1, \ldots, 0) \ominus (n, 0, \ldots, 0) = (0, n-1, n-2, \ldots, 0) \in \widetilde{\operatorname{\textnormal{\L}}_{n+1}},$$
we have
$$(0, n-1, n-2, \ldots, 0) \ominus (n-2, n-2, \ldots, n-2) = (0, 1, 0, \ldots, 0) \in \widetilde{\operatorname{\textnormal{\L}}_{n+1}}.$$

With a similar argument we get
$$
(1, 0, \ldots, 0), (0, 1, 0, \ldots, 0), (0, 0, 1, \ldots, 0), \ldots, (0, \ldots, 0, 1) \in \widetilde{\operatorname{\textnormal{\L}}_{n+1}}.
$$
Since every function $f$ can be written as
$$f = (f(0), 0, \ldots, 0) \oplus (0, f(1), 0, \ldots, 0) \oplus \ldots \oplus (0, 0, \ldots, 0, f(n)),$$
we have $\operatorname{\textnormal{\L}}_{n+1}^{\operatorname{\textnormal{\L}}_{n+1}} \subseteq \widetilde{\operatorname{\textnormal{\L}}_{n+1}}$.
\end{proof}

\begin{proposition}\label{01tilde}
The CMV-algebra $\widetilde{[0,1]}$ consists of the continuous piecewise linear functions from $[0,1]$ to $[0,1]$ where each piece has the form $x \longmapsto mx + \alpha$, with $m \in \mathbb{Z}$ and $\alpha \in \mathbb{R}$.
\end{proposition}
\begin{proof}
Let us denote by $A$ the set of continuous piecewise linear functions of the form above. It is easy to see that $\mathbf{A} = \langle A, \oplus, *, 0, \circ, i \rangle$ with pointwise defined operations and with composition of functions is a CMV-subalgebra of $[0,1]^{[0,1]}$ containing the constant functions. Note that the MV-subalgebra of $[0,1]^{[0,1]}_{\mathcal{MV}}$ generated by constant functions $c_r$, with $r \in [0,1]$, and the identity is precisely $\mathbf{A}_{\mathcal{MV}} = \langle A, \oplus, *, 0, 1\rangle$, hence $\mathbf{A}$ is the smallest CMV-algebra containing $c_r$ for any $r \in [0,1]$.
\end{proof}

\begin{definition}\label{def:constant}
An element $a$ of a CMV-algebra $A$ is a \emph{constant} if $a \diamond x = a$ for all $x \in A$.
\end{definition}
For example $0$ and $1$ are constants. Constant functions in the CMV-algebra $\mathbf{M}^M$ are constants.

\begin{lemma}\emph{\cite{diflge}}
The following properties of constants hold.
\begin{itemize}
\item[(i)] An element $a$ of a CMV-algebra $\mathbf{A}$ is a constant if and only if $a \diamond 1 = a$ if and only if $a \diamond 0 = a$.
\item[(ii)] If $f: \mathbf{A} \longrightarrow \mathbf{B}$ is a morphism of CMV-algebras and $a \in A$ is a constant, then $f(a) \in B$ is a constant. \item[(iii)] The set $K$ of all the constants of the CMV-algebra $\mathbf{A}$ is an MV-subalgebra of $\mathbf{A}_{\mathcal{MV}}$ and an ideal of the monoid $\langle  A, \diamond, i\rangle$.
\end{itemize}
\end{lemma}

\section{CMV-algebras and MV-endomorphisms}\label{cmvend}

If $\mathbf{A} = \langle A, \oplus, ^*, 0\rangle$ is an MV-algebra, then $E(\mathbf{A})$ shall denote the monoid of MV-endomorphisms with composition operation $\boxdot$ given for every $a \in A$ by $(f \boxdot g)(a) = g(f(a))$, and the identity map, $\operatorname{id}_A$, as neutral element. With the above notations we have

\begin{proposition}\label{endcomp}
Let $\mathbf{A} = \langle A, \oplus, ^*, 0, \diamond, i\rangle$ be a CMV-algebra. Then there exists an injective homomorphism of monoids from $\mathbf{A}_\mathcal{M}$ to $E(\mathbf{A}_{\mathcal{MV}})$.
\end{proposition}
\begin{proof}
For all $x, y \in A$ we set $\mu_y(x) = x \diamond y$. It is easy to see that for any $y \in A$, $\mu_y$ is an MV-endomorphism, i.e. $\mu_y \in E(\mathbf{A}_{\mathcal{MV}})$. Then we can define a map $\mu: \mathbf{A}_\mathcal{M} \longrightarrow E(\mathbf{A}_{\mathcal{MV}})$, by setting $\mu(y) = \mu_y$, for all $y \in A$. It is easy to check that $\mu$ is a homomorphism of monoids. To show that $\mu$ is injective we observe that, assuming $\mu(y) = \mu(z)$, we have $x \diamond y = x \diamond z$ for all $x\in A$. So, for $x = i$, we have $y = z$.
\end{proof}

\begin{proposition}\emph{\cite{diflge}}\label{prop226}
Let $\mathbf{A} = \langle A, \oplus, ^*, 0\rangle$ be an MV-algebra and $M_A = \langle A, \diamond, i\rangle$ be any monoid over $A$ such that the map $\Psi: M_A \longrightarrow E(\mathbf{A})$, given by $\Psi(y)(x) = x \diamond y$ for all $x, y \in A$ is a monoid homomorphism. Then we have:
\begin{enumerate}
\item[$(i)$] $\mathbf{A} = \langle A, \oplus, ^*, 0, \diamond, i\rangle$ is a CMV-algebra;
\item[$(ii)$] $\mathbf{A}_{\mathcal{MV}} = \mathbf{A}$;
\item[$(iii)$] $M_A = \mathbf{A}_\mathcal{M}$;
\item[$(iv)$] $\Psi$ is injective.
\end{enumerate}
\end{proposition}

\begin{proposition}\emph{\cite{diflge}}\label{prop:MVtoCMV}
Let $\mathbf{A} = \langle A, \oplus, ^*, 0\rangle$ be an MV-algebra such that there is a mapping $\Phi: A \longrightarrow E(\mathbf{A})$ satisfying the following conditions:
\begin{enumerate}
\item[$(i)$] $\Phi$ is injective;
\item[$(ii)$] $\Phi[A]$ is a submonoid of $E(\mathbf{A})$;
\item[$(iii)$] $\Phi(\Phi(z)(y))(x) = (\Phi(z) \boxdot \Phi(y))(x)$, for all $x, y, z \in A$;
\item[$(iv)$] $\Phi(x)(i) = x$, for any $x\in A$, where $i = \Phi^{-1}(\operatorname{id}_A)$.
\end{enumerate}
Then there exists a CMV-algebra $\mathbf{B}$ such that $\mathbf{A} = \mathbf{B}_{\mathcal{MV}}$.
\end{proposition}

Consider the free MV-algebra $\mathbf{M}_n$ over $n$ generators $p_1, \ldots, p_n$. For any $f \in M_n$ let $\sigma_f: M_n \longrightarrow M_n$ be the map defined by setting $\sigma_f(p_1) = f$ and $\sigma_f(p_i) = p_i$ for $i = 2, \ldots, n$. Then $\sigma_f$ is an endomorphism of $\mathbf{M}_n$ mapping each element $g(p_1, \ldots, p_n)$ of $M_n$ in the element $g(f, p_2, \ldots, p_n)$ of $M_n$ obtained by substituting each occurrence of $p_1$ in $g$ with $f$. We will call the endomorphisms $\sigma_f$ \emph{partial substitutions}.

The map $\Phi: f \in M_n \longmapsto \sigma_f \in E(M_n)$ satisfies the assumptions of Proposition~\ref{prop:MVtoCMV}. Indeed $\Phi$ is clearly an injective function and, for all $f, g \in M_n$, $\sigma_f \circ \sigma_g = \sigma_{\sigma_f(g)}$. Since $\operatorname{id}_{M_n} = \sigma_{p_1}$, then $\Phi[M_n]$ is a submonoid of $E(M_n)$ and $\Phi(\Phi(h)(g))(f) = \Phi(\sigma_h(g))(f) = \sigma_{\sigma_h(g)}(f) = (\sigma_h \circ \sigma_g)(f) (\Phi(h)\circ \Phi(g))(f)$. Further, $\Phi^{-1}(\operatorname{id}_{M_n}) = p_1$ and $\sigma_f(p_1) = f$ by definition. Hence $M_n$ can be equipped with a structure of CMV-algebra by setting $f \diamond g = \sigma_g(f)$ and $i = p_1$.

\section{Congruences, CMV-ideals and simple algebras}\label{cmvcong}

\begin{definition}
An MV-ideal $I$ of the MV-reduct $\mathbf{A}_{\mathcal{MV}}$ of a CMV-algebra $\mathbf{A}$ is a \emph{CMV-ideal} if it is also a right ideal of the monoid $\langle A, \diamond, i\rangle$ and $x \sim_I y$ implies $a \diamond x \sim_I a \diamond y$ for all $a, x, y \in A$, where $\sim_I$ is the MV-congruence generated by $I$ (see Definition~\ref{mvid}).

Alternatively, we can say that a non-empty subset of a CMV-algebra $\mathbf{A}$ is a CMV-ideal if the following conditions are satisfied:
\begin{enumerate}
\item[$(i)$] if $x \in I$ and $x \geq y$ then $y \in I$,
\item[$(ii)$] if $x, y \in I$ then $x \oplus y \in I$,
\item[$(iii)$] $x \diamond y \in I$ for all $x \in I$ and $y \in A$,
\item[$(iv)$] for all $a, x, y \in A$, if $x \sim_I y$ then $a \diamond x \sim_I a \diamond y$.
\end{enumerate}

A subset of $\mathbf{A}$ satisfying conditions $(i$--$iii)$ is called a \emph{$\diamond$-ideal}.
\end{definition}

\begin{lemma}\label{idcong}
Let $\mathbf{A}$ be a CMV-algebra and $I$ a subset of $A$. Then $I$ is a CMV-ideal if and only if the relation $\sim_{I}$ is a CMV-algebra congruence.
\end{lemma}
\begin{proof}
From the definition of CMV-ideal it follows that, if $x \sim_I y$ then $a \diamond x \sim_I a \diamond y$ for all $a \in A$, hence $\sim_I$ is a CMV-congruence.

Conversely, if $\sim_I$ is a congruence of CMV-algebras then, for $x \in I$ and $y \in A$, we have $x \sim_I 0$ and $(x \diamond y) \sim_I (0 \diamond y) = 0$, i.e. $(x \diamond y) \sim_I 0$; hence $(x \diamond y) \in I$. Moreover, to show that $I$ is a CMV-ideal, assume $x \sim_I y$. Since $\sim_I$ is a congruence, we have, for any $a \in A$, $a \diamond x \sim_I a \diamond y$. Hence $I$ is a CMV-ideal.
\end{proof}

The following result has a standard algebraic proof.
\begin{proposition}
Let $\mathbf{A}$, $\mathbf{B}$ be CMV-algebras and $h: A \longrightarrow B$ a CMV-homomorphism. Then for any CMV-ideal $J$ of $\mathbf{B}$ the set $h^{-1}(J) = \{x \in A \ \mid \ h(x) \in J\}$ is a CMV-ideal of $\mathbf{A}$.
\end{proposition}

\begin{exm}\label{ex:3.9}
\begin{enumerate}
\item Consider the CMV-algebra $\mathbf{M}^\mathcal C_1$ of McNaughton functions of one variable. Let $I$ be the set of all the functions $f \in M_1$ such that $f(0) = f(1) = 0$. Then $I$ is a CMV-ideal.

\emph{Indeed, it is clear that $I$ is an MV-ideal; now if $f \in I$ and $g \in M_1$ then $g[\{0,1\}] \subseteq \{0,1\}$ and therefore $f \circ g \in I$. It is easy to see that $f \sim_I g$ iff $f(0) = g(0)$ and $f(1) = g(1)$. Then we have that $h \circ f \sim h \circ g$ for any $h \in M_1$.}
\item To get another example of CMV-ideal on $\mathbf{M}^\mathcal C_1$, consider the MV-ideal $J_0$ of the functions that are constantly equal to zero in a
neighbourhood of $0$ (the neighbourhood depending on the function) and the MV-ideal $J_1$ constructed in the same way for neighbourhoods of $1$. They are prime MV-ideals. Let $J$ be $J_0 \cap J_1$. Then $J$ is a CMV-ideal.

\emph{Clearly $J$ is an MV-ideal; let $f \in J$ and suppose $f = 0$ on $[0,c) \cup (d,1]$. Take $h \in M_1$ and assume, for example, that $h(0) = 0$ and $h(1) = 0$. Since $h$ is continuous, it is possible to find a neighbourhood $U$ of $0$ and a neighbourhood $V$ of $1$ such that $h[U], h[V] \subseteq [0,c)$ and, then, $f \circ h \in J$. The other cases are treated in the same way. Now if $\sim_J$ is the equivalence defined by $J$ we have that $f \sim_J g$ iff $f = g$ in a neighbourhood of $0$ and in a neighbourhood of $1$. It is obvious that, in this case, $h \circ f \sim_J h \circ g$ for any $h \in M_1$.}
\end{enumerate}
\end{exm}

We generalize Example~\ref{ex:3.9} as follows. Let $\mathbf{A}$ be a CMV-algebra and $\mathbf{B}$ an MV-subalgebra of $\mathbf{A}$. We set $S_B = \{a \in A \ \mid \ a \diamond x \in B \ \forall x \in B\}$. It is immediate to verify that $S_B$ is a CMV-subalgebra of $A$. The interpretation in ``functional'' terms is: $S_B$ is the CMV-algebra of the functions that fix the subalgebra $B$. Now we set $J = \{a \in A \ \mid \ a \diamond x = 0 \ \forall x \in B\}$. Obviously $J \subseteq S_B$. Moreover $J$ is a CMV-ideal of $S_B$.

As concrete examples we can consider the CMV-algebra $\operatorname{\textnormal{\L}}_{n+1}^{\operatorname{\textnormal{\L}}_{n+1}}$ and take the subalgebra $B$ of the constant functions. In this case the quotient CMV-algebra $S_B/J$ can be identified with the CMV-subalgebra of $\operatorname{\textnormal{\L}}_{n+1}^{\operatorname{\textnormal{\L}}_{n+1}}$ consisting of the functions $f$ satisfying $f[\{0,1\}] \subseteq \{0,1\}$. It is a CMV-algebra with $4(n+1)^{n-1}$ elements.

According to Universal Algebra, we define \emph{simple} CMV-algebras as those algebras having $\{0\}$ and the whole algebra as the only CMV-ideals.

\begin{lemma}\label{lem:id}
Let $\mathbf{M}$ be an MV-algebra and $\mathbf{A} = \mathbf{M}^M$ the corresponding functional CMV-algebra. Let $I$ be a $\diamond$-ideal of $A$, $a \in M$ and $I(a) = \{f(a) \ \mid \ f \in I\} \subseteq M$. Then $I(a)$ is an MV-ideal of $\mathbf{M}$.
\end{lemma}
\begin{proof}
If $x, y \in I(a)$ then $x = f(a), y = g(a)$ with $f, g \in I$; $f \oplus g \in I$ so $x \oplus y = (f \oplus g)(a) \in I(a)$ and if $z \leq f(a)$ then the function $h: M \longrightarrow M$ defined by $h(a) = z, h(b) = f(b)$ for $b \neq a$ is in $I$ and so $z \in I(a)$. We have more: $I(a) = \bigcup_{f \in I} f[M]$ (hence it is independent from $a$). This is due to the fact that, denoting by $c_b$ the function constantly equal to $b$, if $I$ is a $\diamond$-ideal, then $f \diamond c_b = c_{f(b)} \in I$.
\end{proof}

\begin{theorem}\emph{\cite{diflge}}\label{th:simple}
If $\mathbf{M}$ is an MV-algebra then $\mathbf{A} = \mathbf{M}^M$ is a simple CMV-algebra.
\end{theorem}

Let $\mathbf{A} = \langle A, \oplus, ^*, 0\rangle$ be an MV-algebra and $\mathbf{B}$ a CMV-subalgebra of the CMV-algebra $\mathbf{A}^\mathbf{A}$. The zero element of $B$ will be denoted by $f_{0}$. For any element $f$ of $B$, we set $\operatorname{zero}(f) = \{a \in A \ \mid \ f(a) = 0\}$ and we denote by $\operatorname{supp}(f)$ the complement of $\operatorname{zero}(f)$. For a $\diamond$-ideal $I$ of $\mathbf{B}$ we set $\Sigma(I) = \bigcap_{f \in I} \operatorname{zero}(f)$. Let $S$ be any subset of $A$, then we set $Z(S) = \{f \in B \ \mid \ f(s) = 0 \quad \forall s \in S\}$. $S$ will be called $B$-\emph{stable} if and only if, for every $f \in B$, $f[S] \subseteq S$. With the above notations and definitions we have
\begin{lemma}\label{lem1}
Let $S$ be a $B$-stable non-empty subset of $A$. Then $Z(S)$ is a proper CMV-ideal of $\mathbf{B}$.
\end{lemma}
\begin{proof}
It is easy to check that $Z(S)$ is a proper MV-ideal. Moreover, for $f \in Z(S)$ and $g \in B$ we have, for every $s \in S$, $g(s) \in S$ and $f(g(s)) = 0$. Thus $f \diamond g \in Z(S)$, i.e. $Z(S)$ is a $\diamond$-ideal. Finally, for $f, g \in B$, assume $f \sim_{Z(S)} g$. Then $d(f,g) \in Z(S)$, that is $f(s) = g(s)$ for all $s \in S$. Hence, for any $h \in B$ and $s \in S$, $d(h \diamond f, h \diamond g)(s) = 0$. Therefore $h \diamond f \sim_{Z(S)} h \diamond g$ and $Z(S)$ is a CMV-ideal.
\end{proof}

\begin{lemma}\label{lem2}
Let $J$ be a proper $\diamond$-ideal of $\mathbf{B}$ such that $\Sigma(J) \neq \varnothing$. Then $\Sigma(J)$ is a $B$-stable subset of $A$.
\end{lemma}
\begin{proof}
Let $g\in J$, then for an arbitrary element $h$ of $B$ we have $g \diamond h \in J$ and, for any $a \in \Sigma(J)$, $g(h(a)) = 0$. Hence $h(a) \in \operatorname{zero}(g)$. So $h(a) \in \Sigma(J)$, i.e. $h(\Sigma(J)) \subseteq \Sigma(J)$.
\end{proof}

\section{Modules}\label{mod}

In the present section we are going to give the definition of CMV-module. The main result of the section states that every MV-algebra $\mathbf{B}$ is an
$\mathbf{M}^\mathcal C_1$-module, where $\mathbf{M}^\mathcal C_1$ is the CMV-algebra of one-variable McNaughton functions. In particular, if $\mathbf{B}$ is the MV-algebra $\mathbf{M}_1$, we can associate, to each element $f$ of the MV-algebra $\mathbf{M}_1$, an $\mathbf{M}^\mathcal C_1$-module endomorphism $\Phi_f$ via the mapping $\Phi_f(h) = f \circ h$. In this way the algebraic nature of substitutions is displayed in the suitable framework of $\mathbf{M}^\mathcal C_1$-module.

\begin{definition}\label{def:5.1}
Let us consider an MV-algebra $\mathbf{M} = \langle M, \oplus_M, ^{*_M}, 0_M, 1_M\rangle$ and a CMV-algebra $\mathbf{A} = \langle A, \oplus_A, ^{*_A}, 0_A, \diamond, i\rangle$. $\mathbf{M}$ is an \emph{$\mathbf{A}$-module} if there is an external law
$$\begin{array}{ccc}
A \times M & \longrightarrow & M \\
(a,x) & \longmapsto & ax
\end{array}$$
such that
\begin{enumerate}
\item[$(i)$] $(a \oplus_A b) x = ax \oplus_M bx$,
\item[$(ii)$] $a^{*_A} x = (ax)^{*_M }$,
\item[$(iii)$] $0_A x = 0_M$,
\item[$(iv)$] $ix = x$,
\item[$(v)$] $(a \diamond b) x = a(bx)$,
\end{enumerate}
for all $a, b \in A$ and $x \in M$.
\end{definition}
If $\mathbf{M}$ and $\mathbf{P}$ are $\mathbf{A}$-modules then a module morphism from $\mathbf{M}$ to $\mathbf{P}$ is a morphism of MV-algebras $\varphi: M \longrightarrow P$ such that $\varphi(ax) = a \varphi(x)$ for all $a \in A$, $x \in M$. It is a trivial remark that $\mathbf{A}$-modules and their morphisms form a category in a canonical way.

\begin{exm}
\begin{enumerate}
\item The MV-algebra reduct $\mathbf{A}_{\mathcal{MV}} = \langle A, \oplus, ^*, 0\rangle$ of a CMV-algebra $\mathbf{A}$ is an $\mathbf{A}$-module by setting $a x = a \diamond x$.
\item Let $\mathbf{M}$ be an MV-algebra and $\mathbf{A} = \mathbf{M}^M$. We define $fx = f(x)$ for all $f \in A$ and $x \in M$. This clearly makes $\mathbf{M}$ an $\mathbf{A}$-module.
\item One can generalize the previous example by taking a CMV-subalgebra $\mathbf{B}$ of $\mathbf{A}$. So, for example, $[0,1]$ is a module over the CMV-algebra $\mathbf{M}^\mathcal C_1$ of McNaughton functions.
\item An MV-algebra is an $A_4$-module if and only if it is a Boolean algebra. Indeed it is obvious that the external law can be uniquely defined and satisfies all the conditions needed.
\item Let $\mathbf{A}$ be a CMV-algebra and $\mathbf{K}$ be the MV-algebra of the constants of $A$ (see Definition~\ref{def:constant}), it is easy to prove that $K$ is an MV-subalgebra of $A$. Then $\mathbf{K}$ is, in a natural way, an $\mathbf{A}$-module.

Indeed we define $a k = a \diamond k$ for all $a \in A$ and $k \in K$; we have $a k \diamond x = (a \diamond k) \diamond x = a \diamond (k \diamond x) = a \diamond k = ak$ and the definition is correct. Conditions $(i$---$v)$ of Definition~\ref{def:5.1} are trivially satisfied.
\item If $\mathbf{A}$ is a CMV-subalgebra of $\mathbf{B}$ then $\mathbf{B}$ is in a canonical way an $\mathbf{A}$-module.
\item By combining examples 5 and 6, if $\mathbf{K}_B$ is the algebra of constants of $\mathbf{B}$ then $\mathbf{K}_B$ is an $\mathbf{A}$-module.
\item We can generalize Example 6 to the case of a morphism of CMV-algebras $\Phi: A \longrightarrow B$ and consider $B$ as an $\mathbf{A}$-module defining $ab = \Phi(a) \diamond b$ for all $a \in A$ and $b \in B$. It is easy to check conditions of Definition~\ref{def:5.1}. Even more generally, the morphism $\Phi$ turns naturally every $\mathbf{B}$-module into an $\mathbf{A}$-module; this operation is functorial and may be called ``restriction of
scalars''.
\item Let $\mathbf{M}$ be an $\mathbf{A}$-module, $X$ a given (non-empty) set and consider the set $P$ of all the functions $f: X \longrightarrow M$. $P$ is an MV-algebra in an obvious way (pointwise defined operations). Define $af$, for $a \in A$ and $f\in P$, by $(af)(x) = af(x)$, for all $x \in X$. We obtain a structure of $\mathbf{A}$-module on $P$. In particular the MV-algebra $\mathbf{A}^n$ is an $\mathbf{A}$-module.
\item Let $\mathbf{M}$ be an MV-algebra and $X,P$ as in the previous example; let $\mathbf{A}$ be the CMV-algebra $\mathbf{M}^M$. Then $P$ is an $\mathbf{A}$-module by defining $af = a \circ f$ where $\circ$ is the composition of functions.
\end{enumerate}
\end{exm}

Let $\mathbf{M}^\mathcal C_1$ be the CMV-algebra of McNaughton functions of one
variable where the composition is the usual composition of functions
(see Example~\ref{ex:CMV}.4).

If $f \in M_1$ we denote by $\Phi_f: M_1 \longrightarrow M_1$ the function $\Phi_f(h) = h \circ f$. It is clear that $\Phi_f$ is an endomorphism of the $\mathbf{M}^\mathcal C_1$-module $\mathbf{M}_1$. We call \emph{substitutions} the morphisms of type $\Phi_f$. Note that $\Phi_f(i) = f$.

\begin{lemma}\label{lem:5.3}
$\Phi_{f \circ g} = \Phi_g \circ \Phi _f$.
\end{lemma}
\begin{proof}
$\Phi_{f \circ g}(h) = h \circ (f \circ g) = (h \circ f) \circ g = \Phi_g (\Phi_f(h))$.
\end{proof}

Let $\mathbf{B}$ be any MV-algebra. If $x \in B$ then there is a unique morphism of MV-algebras $\varphi_x: \mathbf{M}_1 \longrightarrow \mathbf{B}$ such that $\varphi_x(i) = x$. This is due to the fact that the MV-algebra $\mathbf{M}_1$ is free over the generator $i$, hence elements of $M_1$ are MV-terms modulo provable equality in $i$ and, by the universal property of free algebras, for any $x \in B$ the map $\varphi_x$ sends any MV-term $\tau(i)$ to the element $\tau(x) \in B$. Denote also $fx = \varphi_x(f)$ for all $f \in M_1$ and $x \in B$.

\begin{proposition}\label{prop5.4}
$\varphi_x \circ \Phi_f = \varphi_{fx}$.
\end{proposition}
\begin{proof}
Elements of the free algebra $M_1$ are MV-polynomials in $i$. By the universal property of free algebras, for any $x \in B$ the map $\varphi_x$ sends any polynomial $\tau(i)$ into the element $\tau(x) \in B$. Hence, for all $g \in M_1$, $(\varphi_x \circ \Phi_f)(g) = \varphi_x(g \circ f) = g(f(x)) = \varphi_{fx}(g)$.
\end{proof}

\begin{theorem}
Every MV-algebra $\mathbf{B}$ is, in a canonical way, a module over the CMV-algebra $\mathbf{M}^\mathcal C_1$ of McNaughton functions.
\end{theorem}
\begin{proof}
We observe that the function $(f,x) \longmapsto fx$ gives to $\mathbf{B}$ a structure of $\mathbf{M}^\mathcal C_1$-module. Then it is enough to prove that $(f \circ g)x = f(gx)$ for all $f, g \in M_1$ and $x \in B$. We have, by Proposition~\ref{prop5.4}, $\varphi_{(f \circ g)x} = \varphi_x \circ \Phi_{f \circ g}$ and, using Lemma~\ref{lem:5.3}, $\varphi_{(f \circ g)x} = \varphi_x \circ (\Phi_g \circ \Phi_f) = (\varphi_x \circ \Phi_g) \circ \Phi_f = \varphi_{gx} \circ \Phi_f = \varphi_{f(gx)}$. Now we can apply all this to $i$ and get the result.
\end{proof}

\section{The logic $\operatorname{S}\operatorname{\textnormal{\L}}_\omega^1$}\label{logic}

We shall extend the fragment of \L ukasiewicz propositional calculus of one variable $v$, here denoted by $\operatorname{\textnormal{\L}}_\omega^1$. The extension of $\operatorname{\textnormal{\L}}_\omega^1$ that we are going to define will be denoted by $\operatorname{S}\operatorname{\textnormal{\L}}_\omega^1$; formulas of  $\operatorname{S}\operatorname{\textnormal{\L}}_\omega^1$ are built from the connectives of negation ($\neg $) and implication $\to$ in the usual way (the conjunction $\odot$ and the disjunction $\oplus$ being considered as derived from $\neg$ and $\to$), and by a binary connective $\blacktriangleleft$ defined as follows:

$\varphi \blacktriangleleft \psi $ is the formula obtained by the formula $\varphi$ where every occurrence of the variable $v$ is substituted by the formula $\psi$. Actually the set $\operatorname{Form}$ of formulas of $\operatorname{S}\operatorname{\textnormal{\L}}_\omega^1$ is obtained inductively by the following stipulations:
\begin{enumerate}
\item The single variable $v$ is a (atomic) formula;
\item If $\alpha$ is a formula, then $\neg \alpha$ is a formula;
\item If $\alpha$ and $\beta$ are formulas, then $\alpha \rightarrow \beta$ and $\alpha \blacktriangleleft \beta$ are formulas.
\end{enumerate}

In the light of next results, the formula $\alpha \blacktriangleleft \beta$ can be thought of as the formula obtained from $\alpha$ by the substitution of each occurrence of the variable $v$ with the formula $\beta$.

Interpretation of connectives of $\operatorname{S}\operatorname{\textnormal{\L}}_\omega^1$ is given by the following definition.

\begin{definition}
Let $\mathbf{A} = \langle A, \oplus, ^*, 0, \diamond, i \rangle$ be a CMV-algebra. Then an \emph{$\mathbf{A}$-valuation} is a function $\operatorname{val}: \operatorname{Form} \longrightarrow A$ satisfying the following properties, where $\alpha$ and $\beta$ denote arbitrary formulas:
\begin{itemize}
\item[-] $\operatorname{val}(v) = i$,
\item[-] $\operatorname{val}(\neg \alpha) = (\operatorname{val}(\alpha))^*$,
\item[-] $\operatorname{val}(\alpha \rightarrow \beta) = \operatorname{val}(\alpha) \rightarrow \operatorname{val}(\beta))$,
\item[-] $\operatorname{val}(\alpha \blacktriangleleft \beta) = \operatorname{val}(\alpha) \diamond \operatorname{val}(\beta)$.
\end{itemize}
\end{definition}

An $\mathbf{A}$-valuation $\operatorname{val}$ is said to $\mathbf{A}$-satisfy a formula $\alpha$ iff $\operatorname{val}(\alpha) = 1$; $\alpha$ is an $\mathbf{A}$-tautology iff $\alpha$ is satisfied by all $\mathbf{A}$-valuations. Formulas $\alpha$ and $\beta$ are semantically $\mathbf{A}$-equivalent iff $\operatorname{val}(\alpha) = \operatorname{val}(\beta)$ for all $\mathbf{A}$-valuations $\operatorname{val}$. The syntactic approach to many-valued logics is the same as for propositional classical logic. A set of formulas, called {\em axioms}, and a set of inference rules are fixed. In $\operatorname{S}\operatorname{\textnormal{\L}}_\omega^1$ we have

Axioms:
\begin{enumerate}
\item[$\operatorname{Ax}1$] \quad $\varphi \to (\psi \to \varphi)$
\item[$\operatorname{Ax}2$] \quad $(\varphi \to \psi) \to ((\psi \to \chi)\to (\varphi\to \chi))$
\item[$\operatorname{Ax}3$] \quad $((\varphi \to \psi)\to \psi)\to ((\psi \to \varphi)\to \varphi)$
\item[$\operatorname{Ax}4$] \quad $(\neg \varphi\to \neg \psi)\to (\psi \to \varphi)$
\item[$\operatorname{Ax}5$] \quad $(\varphi \blacktriangleleft v)\to \varphi$
\item[$\operatorname{Ax}6$] \quad $(v\blacktriangleleft \varphi)\to \varphi$
\item[$\operatorname{Ax}7$] \quad $\varphi \to (v \blacktriangleleft \varphi)$
\item[$\operatorname{Ax}8$] \quad $\varphi \to (\varphi \blacktriangleleft v)$
\item[$\operatorname{Ax}9$] \quad $((\varphi \to \psi)\blacktriangleleft \gamma) \to ((\varphi\blacktriangleleft \gamma) \to (\psi \blacktriangleleft \gamma))$
\item[$\operatorname{Ax}10$] \quad $((\neg \varphi \to \psi)\blacktriangleleft \gamma)\to (\neg(\varphi\blacktriangleleft \gamma) \to (\psi \blacktriangleleft \gamma))$
\item[$\operatorname{Ax}11$] \quad $(\neg(\varphi\blacktriangleleft \gamma) \to (\psi \blacktriangleleft \gamma))\to ((\neg \varphi \to \psi)\blacktriangleleft \gamma)$
\item[$\operatorname{Ax}12$] \quad $(\varphi \blacktriangleleft (\psi\blacktriangleleft \gamma))\to ((\varphi\blacktriangleleft \psi)\blacktriangleleft \gamma)$
\item[$\operatorname{Ax}13$] \quad $((\varphi\blacktriangleleft \psi)\blacktriangleleft \gamma)\to (\varphi \blacktriangleleft (\psi\blacktriangleleft \gamma))$
\item[$\operatorname{Ax}14$] \quad $\neg(\varphi\blacktriangleleft \psi)\to (\neg\varphi\blacktriangleleft \psi)$
\item[$\operatorname{Ax}15$] \quad $(\neg\varphi\blacktriangleleft \psi)\to \neg(\varphi\blacktriangleleft \psi)$
\end{enumerate}

Inference rules:
$$\frac{\alpha \quad \alpha \to \beta}{\beta} \qquad \textrm{\emph{Modus Ponens}},$$
$$\frac{\alpha}{\alpha \blacktriangleleft \beta} \qquad \textrm{$\blacktriangleleft$-\emph{rule}},$$
$$\frac{\alpha}{\alpha \to (\alpha \blacktriangleleft \beta)} \qquad \textrm{\emph{Arrow $\blacktriangleleft$-rule}}.$$
By $\vdash \alpha $ we mean that $\alpha$ is a theorem of $\operatorname{S}\operatorname{\textnormal{\L}}_\omega^1$, and the set of all the theorems of $\operatorname{S}\operatorname{\textnormal{\L}}_\omega^1$ shall be denoted by $\operatorname{Th}$. Let now $\sim_{\blacktriangleleft}$ be the binary relation on $\operatorname{Form}$ defined by:
$$\alpha \sim_{\blacktriangleleft} \beta \qquad \text{iff} \qquad \vdash \alpha \to \beta \ \text{ and } \ \vdash \beta \to \alpha.$$

With the above notations we have:
\begin{proposition}\label{simequiv}
$\sim_{\blacktriangleleft}$ is an equivalence relation such that if $\alpha \sim_{\blacktriangleleft} \beta$ then $((\alpha \blacktriangleleft \gamma) \sim_{\blacktriangleleft} (\beta \blacktriangleleft \gamma))$, for every formula $\gamma$.
\end{proposition}
\begin{proof}
$\sim_{\blacktriangleleft}$ is trivially an equivalence relation. The fact that $\alpha \sim_{\blacktriangleleft} \beta$ implies $((\alpha \blacktriangleleft \gamma) \sim_{\blacktriangleleft} (\beta \blacktriangleleft \gamma))$, for all $\gamma$, is as an easy application of $\blacktriangleleft$-rule, $\operatorname{Ax}9$ and Modus Ponens.
\end{proof}

The equivalence classes of $\sim_{\blacktriangleleft}$ will be denoted by
$$\mid \alpha \mid \ := \ \{\beta \in \operatorname{Form} \mid \beta \sim_{\blacktriangleleft}\alpha\}.$$

\begin{lemma}\label{thms}
For each formula $\alpha$, $\vdash \alpha$ if and only if $\mid \alpha \mid \ = \ \operatorname{Th}$.
\end{lemma}
\begin{proof}
If $\mid \alpha \mid = \operatorname{Th} $, then $\vdash \alpha$, because $\alpha \in \mid \alpha \mid$.

Conversely, suppose that $\vdash \alpha$ and $\vdash \beta$. Then, by $\operatorname{Ax}1$ and Modus Ponens, we get $\vdash \beta \to \alpha$ and $\vdash \alpha \to \beta$, hence $\beta \sim_{\blacktriangleleft} \alpha$. Then $\operatorname{Th} \subseteq \mid\alpha\mid$. Now take $\beta \in \mid\alpha\mid$. Then $\vdash \alpha \to \beta$ and therefore, via Modus Ponens, we get $\vdash \beta$, which implies $\mid \alpha \mid \subseteq \operatorname{Th}$.
\end{proof}

\begin{theorem} The quotient set $\frac{\operatorname{Form}}{\sim_{\blacktriangleleft}}$ becomes a CMV-algebra,
once equipped with the operations $\to, \neg, \diamond$ and the
constants $1$ and $i$ given by the following stipulations:
\begin{enumerate}
\item[$(i)$] $\mid \alpha \mid \to \mid \beta \mid \ := \ \mid \alpha \to \beta \mid$
\item[$(ii)$] $(\mid \alpha \mid )^* \ := \ \mid \neg \alpha \mid$
\item[$(iii)$] $\mid \alpha \mid \diamond \mid \beta \mid \ := \ \mid \alpha \blacktriangleleft \beta \mid$
\item[$(iv)$] $1 \ := \ \operatorname{Th}$
\item[$(v)$] $i \ := \ \mid v \mid$.
\end{enumerate}
\end{theorem}
\begin{proof}
By Proposition~\ref{simequiv} and Lemma~\ref{thms} $(i$--$iii)$ are well defined operations on $\frac{\operatorname{Form}}{\sim_{\blacktriangleleft}}$. It remains to prove that the axioms of CMV-algebras are satisfied. Here we limit ourselves to check the proper axioms of CMV-algebra. Let us prove the associativity of $\diamond$. Indeed, to prove that $(\mid \alpha \mid \diamond \mid \beta \mid)\diamond \mid \gamma \mid = \mid \alpha \mid \diamond(\mid \beta \mid \diamond \mid \gamma \mid)$ we have to prove that $\alpha \blacktriangleleft (\beta \blacktriangleleft \gamma) \sim_{\blacktriangleleft} (\alpha \blacktriangleleft \beta) \blacktriangleleft \gamma$. Equivalently we should show that
\begin{equation}\label{(j)}
\vdash (\alpha \blacktriangleleft (\beta \blacktriangleleft \gamma)) \to ((\alpha \blacktriangleleft \beta) \blacktriangleleft \gamma),
\end{equation}
\begin{equation}\label{(jj)}
\vdash ((\alpha \blacktriangleleft \beta) \blacktriangleleft \gamma) \to (\alpha \blacktriangleleft (\beta \blacktriangleleft \gamma)),
\end{equation}
but (\ref{(j)}) is precisely $\operatorname{Ax}12$ and (\ref{(jj)}) is $\operatorname{Ax}13$. To prove the distributivity of $\diamond$ with respect to the $\oplus$ operation we need to prove that
$$(\neg \alpha \to \beta) \blacktriangleleft \gamma \quad \sim_{\blacktriangleleft} \quad \neg (\alpha \blacktriangleleft \gamma) \to (\beta \blacktriangleleft \gamma),$$
i.e. we need to show that the following hold:
\begin{equation}\label{(jjj)}
\vdash ((\neg \alpha \to \beta) \blacktriangleleft \gamma) \to (\neg(\alpha \blacktriangleleft \gamma) \to (\beta \blacktriangleleft \gamma)),
\end{equation}
\begin{equation}\label{(jv)}
\vdash (\neg(\alpha \blacktriangleleft \gamma) \to (\beta \blacktriangleleft \gamma)) \to ((\neg \alpha \to \beta) \blacktriangleleft \gamma).
\end{equation}
Again, (\ref{(jjj)}) and (\ref{(jv)}) are, respectively, $\operatorname{Ax}10$ and $\operatorname{Ax}11$.

To prove that $1 \diamond \mid \alpha \mid = 1$ we have to show that for some $\beta \in \operatorname{Th}$ we have $\beta \sim_{\blacktriangleleft} (\beta \blacktriangleleft \alpha)$, which is absolutely trivial, as well as the fact that $i \diamond \mid \alpha \mid = \mid \alpha \mid$ and $\mid \alpha \mid \diamond i = \mid \alpha \mid$.

Last, to verify that $\mid \alpha \mid^* \diamond \mid \beta \mid = (\mid \alpha \mid \diamond \mid \beta \mid)^*$, we just need to apply $\operatorname{Ax}14$ and $\operatorname{Ax}15$.
\end{proof}

\section{Conclusion}

In this last section we would like to stress some facts about CMV-algebras and $\operatorname{S}\operatorname{\textnormal{\L}}_\omega^1$.

We presented plenty of examples of CMV-algebras and we have seen that the one generated free MV-algebra is in a canonical way a CMV-algebra. Moreover, CMV-algebra have many interesting algebraic properties and connections with MV-algebras. For example, Propositions~\ref{prop226} and~\ref{prop:MVtoCMV} tell us how to extend an MV-algebra $\mathbf{A}$ to a CMV-algebra $\mathbf{A}^\mathcal C$ using the endomorphism monoid of $\mathbf{A}$, and Theorem~\ref{th:simple} highlights a wide class of simple CMV-algebras. 

On the other hand, the system $\operatorname{S}\operatorname{\textnormal{\L}}_\omega^1$, introduced in Section~\ref{logic}, looks like the single variable fragment of some interesting~---~yet unknown~---~expansion of \L ukasiewicz logic. Hence, it could be fruitful to find a suitable $\operatorname{S}\operatorname{\textnormal{\L}}_\omega$ system (i.e. a system in a numerable set of variables) containing $\operatorname{S}\operatorname{\textnormal{\L}}_\omega^1$ as a fragment.

So the theory of CMV-algebras is still at its initial stage, but many signs indicate that it is worth to keep on investigating such structures, their connection with logic and their possible applications (e.g. at the theory of chaotic dynamical systems, as pointed out in~\cite{diflge}).

\end{document}